\input ./style/arxiv-general.cfg
\documentclass[aos,MSNbibl,nameyear,dvips]{arximspdf}
\makeatletter
   \@ifpackageloaded{graphicx}{}{\usepackage{graphicx}}
\makeatother

\doi{10.1214/15-AOS1270E}
\volume{43}
\issue{4}
\pubyear{2015}
\firstpage{1455}
\lastpage{1462}
\docsubty{FLA}
\referstodoi{10.1214/14-AOS1270}

\makeatletter
\newcommand{\rrVert}{\Vert}
\newcommand{\llVert}{\Vert}
\newtheorem{theorem}{Theorem}

\theoremstyle{remark}

\makeatother

\begin{document}
\begin{frontmatter}
\vspace*{12pt}
\title{Discussion of ``Frequentist coverage of adaptive nonparametric
Bayesian credible sets''}
\runtitle{Discussion}

\begin{aug}
\author[A]{\fnms{Subhashis}~\snm{Ghosal}\corref{}\ead[label=e1]{subhashis\_ghoshal@ncsu.edu}}
\runauthor{S. Ghosal}
\affiliation{North Carolina State University}
\address[A]{Department of Statistics\\
North Carolina State University\\
4276 SAS Hall, 2311 Stinson Drive\\
Raleigh, North Carolina 27695-8203\\
USA\\
\printead{e1}}
\end{aug}

%
\received{\smonth{1} \syear{2015}}


\end{frontmatter}

First I would like to congratulate the authors Botond Szab\'o, Aad van
der Vaart and Harry van Zanten for a fine piece of work on the
extremely important topic of frequentist coverage of adaptive
nonparametric credible sets. Credible sets are used by Bayesians to
quantify uncertainty of estimation, which is typically viewed as more
informative than point estimation. Such sets are often easily
constructed, for instance, by sampling from the posterior, while
confidence sets in the frequentist setting may need evaluating limiting
distributions, or resampling, which needs additional justification.
Bayesian uncertainty quantification in parametric problems from the
frequentist view is justified through the Bernstein--von Mises theorem.
In recent years, such results have also been obtained for the
parametric part in certain semiparametric models, guaranteeing coverage
of Bayesian credible sets for it. However, as mentioned by the authors,
inadequate coverage of nonparametric credible sets has been observed
[\citet{Cox1993}, \citet{Freedman1999}] in the white noise model, arguably the
simplest nonparametric model. A clearer picture emerged after the work
of \citet{Knapik2011} that undersmoothing
priors can resolve the issue of coverage; see also \citet{Leahu2011} and
\citet{Castillo2013}.

In the present paper, the authors address the issue of coverage of
credible sets in a white noise model under the inverse problem setting,
when the underlying smoothness (i.e., regularity) of the true parameter
is not known, so a procedure must adapt to the smoothness. The authors
follow an empirical Bayes approach where a key regularity parameter in
the prior is estimated from its marginal likelihood function. As the
authors mentioned, undersmoothing leads to inferior point estimation
and is also difficult to implement when the smoothness of the parameter
is not known. We shall see that the issue of coverage can also be
addressed by two other alternative approaches.

Before entering a discussion on the contents of the paper, let us take
another look at the coverage problem for Bayesian credible sets in an
abstract setting.
Suppose that we have a family of experiments based on observations
$Y^{(n)}$ and indexed by a parameter $\theta\in\Theta$, some
appropriate metric space. Let $\epsilon_n$ be the minimax convergence
rate for estimating $\theta$. Let $\gamma_n\in[0,1]$ be a sequence
which can be fixed or may tend to $0$. For some $m_n\rightarrow\infty$,
typically a slowly varying sequence, the goal is to find a subset
$\mathcal{C}(Y^{(n)})\subseteq\Theta$ such that uniformly on $\theta
_0\in\mathcal{B}$:
\begin{longlist}[(iii)]
\item[(i)] $\Pi(\theta\in\mathcal{C}(Y^{(n)})|Y^{(n)})\ge1-\gamma_n$,
\item[(ii)] $P_{\theta_0}^{(n)}(\theta_0\in\mathcal
{C}(Y^{(n)}))\rightarrow1$,
\item[(iii)] $\operatorname{diam}(\mathcal{C}(Y^{(n)}))=O_{P_{\theta
_0}^{(n)}}(m_n\epsilon_n)$,
\end{longlist}
where $\mathcal{B}$ varies over a class of compact balls in $\Theta$.

In the formulation, credibility may increase with the sample size. We
find it natural that when the information content is increasing, a
researcher should quantify uncertainty with more and more confidence,
instead of staying at a fixed level, just like one seeks for more
precise point estimators or tests. If $\gamma_n\to0$, it can be seen
that the problems of mismatch of credibility and coverage pointed out
in \citet{Cox1993} and \citet{Freedman1999} go away. Thus, although the
uncertainty quantification of a Bayesian and a frequentist may not
match at finite levels, they do match at the infinitesimal level. For
finer matching, one may also like to impose some requirement on how
fast $P_{\theta_0}^{(n)}(\theta_0\in\mathcal{C}(Y^{(n)}))$ should
approach $1$, but we shall forgo the issue in this discussion.
Another approach is to obtain a $(1-\gamma_n)$-credible ball around the
posterior mean typically with fixed $\gamma_n$ and inflate the region
by a factor $m_n$, to be called the inflation factor, to ensure
adequate frequentist coverage. The size of the original credible region
is typically of the order of the minimax convergence rate $\epsilon_n$
so that the third condition will be met.
The factor $m_n$ can be considered as a reasonable price for the
increased level of coverage. Typically, the resulting extra cost $m_n$
is low, for instance, in an asymptotic normality setting, while
adopting $(1-\gamma_n)$-credible sets with $\gamma_n\to0$, the
additional cost is $m_n=o(\sqrt{\log(1/\gamma_n}))$. In the setting we
shall discuss, the inflation factor may be taken as a sufficiently
large constant. The supremum over compact sets in the formulation
imposes honesty of the coverage.

As mentioned by the authors, fully adaptive honest nonparametric
confidence regions are not possible by any means, so in the adaptive
context $\Theta$ will be replaced by an appropriate subset of the
parameter space, such as the set of self-similar sequences or polished
tailed sequences in the context of the paper. The concept of polished
tail is pretty elegant as it blends nicely in the adaptive setting
without any direct reference to the smoothness of the parameter.

The main result proved in the paper, namely, honest coverage of
adaptive posterior credible regions for $\bolds{\theta}=(\theta
_1,\theta
_2,\ldots)$ in the model $Y_i=\kappa_i \theta_i+n^{-1/2} \varepsilon
_i$, where $\bolds{\theta}\in\ell_2$ and $\varepsilon_i\stackrel{\mathrm{i.i.d.}}{\sim}N(0,1)$, for
all polished tail sequences is certainly exciting. In terms of the
equivalent (and perhaps more directly relevant) white noise inverse
problem model $dY(t)=Kf(t)\,dt+n^{-1/2} \,dW(t)$, this translates into
honest coverage of credible regions for $f$ through Parseval's
identity, where the distance on $f$ is measured in terms of the
$L_2$-distance. However, $L_2$-regions for functions do not look like
bands, and may be a little harder to visualize. This aspect may be
relevant for covering a true function that has a bump like the one
given by equation~(4.1) in the paper under discussion, since
$L_2$-closeness does not even imply pointwise closeness, let alone
uniform closeness. Regions on function spaces given by $L_\infty$-neighborhoods are easier to visualize and interpret. Moreover, such
uniform closeness has other implications. For instance, if derivatives
of the functions in a region are uniformly $\epsilon_n$-close to the
derivative of the true function and the derivative of the true function
has a well separated mode, then the mode of a function in that region
is $O(\epsilon_n)$-close to the true mode. The observation can be used
to induce honest confidence regions for the mode from those for the
derivative function under the $L_\infty$-distance.

Study of coverage of $L_\infty$-regions with chosen credibility needs
studying posterior contraction rates under the $L_\infty$-norm, which is
easier if conjugacy is present, like in the white noise model or
nonparametric regression using a random series with normal
coefficients. Below we shall argue that in the white noise model
credible regions for the $L_\infty$-norm can also be characterized and
computed relatively easily, and their coverage can be shown to be
adequate. Interestingly, we can use a fixed level of credibility (any
value higher than $1/2$ works) and the inflation factor can be taken to
be a constant. We shall follow techniques similar to those used in \citet{Yoo2014}, who considered the problem of multivariate
nonparametric regression using a random series of tensor product
B-splines in the known smoothness setting. In a sense the present
treatment of the simpler white noise model will be easier, but there
are certain differences as well, particularly since the number of basis
elements used in constructing the prior is infinite in the present
case, unlike the case treated by \citet{Yoo2014}. For the sake of
simplicity of the discussion, we focus on the direct problem, that
is, $\kappa_i\equiv1$, and consider a Fourier basis $\phi_1(x)=1$,
$\phi_{2i}(x)=\sqrt2\cos(2\pi ix)$, $\phi_{2i+1}(x)=\sqrt2\sin
(2\pi
ix)$, $i=1,2,\ldots.$ Let the true function be denoted by $f_0$ and the
true sequence by $\bolds{\theta}_0=(\theta_{01},\theta_{02},\ldots)$.
Since we intend to study the $L_\infty$-contraction rate and coverage of
$L_\infty$-regions, we need to assume that the true function $f_0$ belongs
to a H\"older class or, stated in terms of coefficients $\sum_{i=1}^{\infty}
i^{\alpha}|\theta_{0i}|<\infty$, which\vspace*{1pt} is stronger than the analogous
Sobolev condition $\sum_{i=1}^{\infty}i^{2\alpha}\theta
_{0i}^2<\infty$.
The logarithmic factor we obtain in the rate is not optimal---it is
off by the factor $(\log n)^{1/2(2\alpha+1)}$. Using a more refined
analysis or perhaps using a different basis like B-splines or wavelets,
an optimal logarithmic factor may be obtained as in \citet{Yoo2014} or \citet{Gine2011}. Also, because we use a Fourier
basis, we assume that the true function is periodic, but this does not
dampen the essential spirit of the argument.

Consider the white noise model $dY(t)=f(t)\,dt+n^{-1/2}\,dW(t)$ and its
equivalent normal sequence model $Y_i=\theta_i+n^{-1/2}\varepsilon_i$,
where $Y_i=\int\phi_i(x)\,dY(x)$, $\theta_i=\int\phi_i(x) f(x)\,dx$ and
$\varepsilon_i=\int\phi_i(x)\,dW(x)$, $i=1,2,\ldots.$
Let the prior $\Pi$ be defined by $\theta_i\stackrel{\mathrm{ind}}{\sim}N(0,i^{-2\alpha+1})$.
Let $\hat f=\mathrm{E}(f|D_n)$, where $D_n$ stands\vspace*{1.5pt} for the data. Note that
$\hat f(x)=\sum_{i=1}^{\infty}\hat{\theta}_i \phi_i(x)$,
where $\hat\theta_i=\mathrm{E}(\theta_i|Y_i)=nY_i/(i^{2\alpha+1}+n)$,
and $\operatorname{var}
(\theta_i|Y_i)=(i^{2\alpha+1}+n)^{-1}$.
Let $B(\alpha,R)=\{f\dvtx  \sum_{i=1}^{\infty}i^{\alpha}|\theta_i|\le
R\}$.
Below we shall write ``$\lesssim$'' for inequality up to a constant and
``$\asymp$'' for equality in order.

\begin{theorem}
\label{main_theorem}
For any $M_n\to\infty$, $\mathrm{E}_{f_0}\Pi_\alpha(f\dvtx  \|f-f_0\|
_\infty>M_n \epsilon
_n|D_n )\to0$
uniformly for all $f_0\in B(\alpha,R)$, where $\epsilon_n =
n^{-\alpha
/(2\alpha+1)}\sqrt{\log n}$ and for a sufficiently large constant
$M>0$, $P_{f_0}\{\|f_0-\hat f\|_\infty\le M h_n \}\to1$ where $h_n$ is
determined by $\Pi_\alpha(f\dvtx  \|f-\hat f\|_\infty\le
h_n|D_n)=1-\gamma$,
$\gamma\ge1/2$ is a predetermined constant. Moreover, $h_n\lesssim
\epsilon_n$.
\end{theorem}

The theorem implies that the $(1-\gamma)$-credible region for
$L_\infty
$-distance around the posterior mean for any $\gamma\le1/2$ inflated
by a sufficiently large factor $M$ has asymptotic coverage $1$ and its
size $h_n$ is not larger than the posterior contraction rate, which is
nearly optimal. It is interesting to note that $h_n$ is actually
deterministic since the posterior distribution of $f-\hat f$ is free of
the observations. Analytical computation of $h_n$ may be difficult, but
can be easily determined by simulations.

\begin{pf*}{Proof of  Theorem~\protect\ref{main_theorem}}
We have $f(x)=\sum_{i=1}^{\infty}\theta_i \phi_i(x)$. Thus, given
$D_n$, $Z=f-\hat
f$ is a mean-zero Gaussian process with covariance kernel
\[
\operatorname{cov} \Biggl(\sum_{i=1}^{\infty}
\theta_i \phi_i(s),\sum_{i=1}^{\infty}
\theta_i \phi_i(t) \Big|D_n \Biggr) = \sum
_{i=1}^{\infty} \bigl(i^{2\alpha+1}+n
\bigr)^{-1}\phi_i(s)\phi_i(t)
\]
and
\begin{eqnarray*}
\mathrm{E} \bigl(\bigl|Z(s)-Z(t)\bigr|^2|D_n \bigr) &=& \sum
_{i=1}^{\infty}\operatorname{var}(
\theta_i|D_n)\bigl| \phi_i(s)-\phi
_i(t)\bigr|^2
\\
&\lesssim& \sum_{i=1}^{\infty}
\bigl(i^{2\alpha+1}+n \bigr)^{-1}i^2 |s-t|^2
\\
&\lesssim& n^{2(\alpha-1)/(2\alpha+1)}|s-t|^2
\end{eqnarray*}
by standard estimates and the fact $|\phi_i(s)-\phi_i(t)|\le2\sqrt
{2}\pi i|s-t|$, a consequence of the mean value theorem and the
boundedness of trigonometric functions. Using a uniform grid with
mesh-width $\delta_n\asymp n^{-p}$ for $p>0$ sufficiently large and a
chaining argument for Gaussian processes with values of $Z$ at the
chosen grid points, Lemma~2.2.2 and Corollary~2.2.8 of \citet{vanderVaart1996} give the estimate $\mathrm{E}\|Z\|_\infty\le\sqrt
{\mathrm{E}\|Z\|_\infty
^2} \lesssim n^{-\alpha/(2\alpha+1)}\sqrt{\log n}$.

Let $V(x)=\hat f(x)-\mathrm{E}_{f_0} \hat f(x)=\sum_{i=1}^{\infty}\sqrt{n}
\varepsilon_i \phi_i(x) /(i^{2\alpha+1}+n)$. Then $V$ is a mean-zero
Gaussian process with covariance kernel $\sum_{i=1}^{\infty
}n(i^{2\alpha+1}+n)^{-2}\times\break 
\phi_i(s)\phi_i(t)$ and
\[
\mathrm{E}\bigl|V(s)-V(t)\bigr|^2 = \sum_{i=1}^{\infty}n
\bigl(i^{2\alpha+1}+n \bigr)^{-2} \bigl|\phi_i(s)-
\phi_i(t)\bigr|^2.
\]
Arguing as before, it follows that $\mathrm{E}_{f_0}\|V\|_\infty\lesssim
n^{-\alpha/(2\alpha+1)}\sqrt{\log n}$.

Now using the uniform boundedness of the basis functions and\break $\sum_{i=1}^{\infty}
i^\alpha|\theta_{0i}|\le R$, uniformly for $f_0\in B(\alpha,R)$, we
have for any $k$, $\sum_{i>k} |\theta_{0i}| \le Rk^{-\alpha}$. Therefore,
\begin{eqnarray*}
\|\mathrm{E}_{f_0}\hat f-f_0\|_\infty&=& \Biggl
\llVert \sum_{i=1}^{\infty
} \biggl(
\frac{n}{i^{2\alpha +1}+n}-1 \biggr)\theta_{0i} \phi _i \Biggr
\rrVert _\infty
\\
& \le& \sqrt2 R \biggl(\frac{k^{\alpha+1}}{n} +k^{-\alpha} \biggr) \lesssim
n^{-\alpha/(2\alpha+1)}
\end{eqnarray*}
by choosing $k=k_\alpha\asymp n^{1/(2\alpha+1)}$.

Combining the three pieces, it follows using Chebyshev's inequality
that the posterior contraction rate under the $L_\infty$-distance is
$\epsilon_n$.

Now we find a lower bound for the size of the credible region.
By definition $h_n$, the $(1-\gamma)$-quantile of the distribution of
$\|Z\|_\infty$ for the mean-zero Gaussian process $Z$ with covariance
kernel $\sum_{i=1}^{\infty}(i^{2\alpha+1}+n)^{-1}\phi_i(s)\phi
_i(t)$ is at least as
large as the median of the distribution of $\|Z\|_\infty$. Now $\sigma
_Z^2=\sup\mathrm{E}|Z(t)|^2$ is easily seen to be $O(n^{-2\alpha
/(2\alpha
+1)})$. Since $\mathrm{E}\|Z\|_\infty^2 \ge\sigma_Z^2$, standard facts about
Gaussian processes imply that $\mathrm{E}\|Z\|_\infty$ and the median of
$\|Z\|_\infty
$ are of the same order [cf. \citet{Ledouz1991}, pages 52 and
54]. Hence, to find a lower bound for $h_n$, it suffices to lower bound
$\mathrm{E}\|Z\|_\infty$. We shall show that the order of the lower bound is
$n^{-\alpha/(2\alpha+1)}\sqrt{\log n}$.

To this end, we observe that $\|Z\|_\infty\ge\max\{Z(j/k_\alpha
):j=1,\ldots,k_\alpha\}$, and
\[
\mathrm{E}\bigl|Z(j/k_\alpha)-Z(l/k_\alpha)\bigr|^2 = \sum
_{i=1}^{\infty
} \bigl(i^{2\alpha+1}+n
\bigr)^{-1} \bigl|\phi _i(j/k_\alpha)-
\phi_i(l/k_\alpha)\bigr|^2.
\]
With a sufficiently small fixed $\epsilon>0$, there exists a $\delta>0$
such that $|\sin s-\sin t|>\epsilon$ if $|s-t|>\delta$ and $|s+t-\pi
|>\delta$, and a similar assertion holds for the cosine function.
Therefore, it is observed that for $j,l=1,\ldots,k_\alpha$, $j\ne l$,
$\phi_i(j/k_\alpha)$ and $\phi_i(l/k_\alpha)$ differ by at least a
fixed positive number for a positive fraction of $i\in\{2,\ldots
,k_\alpha\}$. From this we obtain that there exists $c>0$ such that
\[
\mathrm{E}\bigl|Z(j/k_\alpha)-Z(l/k_\alpha)\bigr|^2\ge c
n^{-2\alpha/(2\alpha+1)}.
\]
Let $U_j=\sqrt2 c^{-1/2} n^{\alpha/(2\alpha+1)}Z(j/k_\alpha)$,
$j=1,\ldots,k_\alpha$, so that $\mathrm{E}(U_j-U_l)^2\ge\mathrm{E}(V_j-V_l)^2$, where
$V_1,\ldots,V_{k_\alpha}\stackrel{\mathrm{i.i.d.}}{\sim}N(0,1)$. Hence, by Slepian's inequality
[cf. Corollary~3.14 of \citet{Ledouz1991}] and equation (3.14)
of \citet{Ledouz1991}, we obtain
\[
\mathrm{E} \Bigl(\max_j U_j \Bigr)\ge
\mathrm{E} \Bigl(\max_j V_j \Bigr)\gtrsim\sqrt{
\log k_\alpha}\asymp\sqrt {\log n},
\]
which upon rescaling gives $\mathrm{E}\|Z\|_\infty\ge\mathrm{E}(\max
Z(j/k_\alpha)
)\gtrsim n^{-\alpha/(2\alpha+1)}\sqrt{\log n}$.

Now turning to coverage,
the lack of coverage of the credible set inflated by a sufficiently
large constant $M$ is given by
\begin{eqnarray*}
P_{f_0}\bigl\{\|f_0-\hat f\|_\infty> M h_n
\bigr\} &\le& \mathrm{P} \bigl\{ \|V\|> M'\epsilon_n -\|
\mathrm{E}\hat f-f_0\|_\infty \bigr\}
\\
& \le& 2 e^{-C \log n}\to0
\end{eqnarray*}
by virtue of Borell's inequality [cf. second assertion of
Proposition~A.2.1 of \citet{vanderVaart1996}],
since $\sup_t \operatorname{var}(Z(t))\lesssim n^{-2\alpha/(2\alpha+1)}$ and
$\epsilon
_n\asymp n^{-\alpha/(2\alpha+1)}\sqrt{\log n}$,
uniformly for $f_0\in B(\alpha,R)$, where $M'$ and $C$ are positive constants.

Finally, we estimate the size of the inflated credible region. For that
we need to find an upper bound for the $(1-\gamma)$-quantile of the
distribution of $\|Z\|_\infty$ given $D_n$. By Borell's inequality [cf.
third assertion of Proposition~A.2.1 of \citet{vanderVaart1996}], it is clear that the $(1-\gamma)$-quantile is bounded by
$\sqrt
{8 \mathrm{E}\|Z\|_\infty^2 \log(2/\gamma)}$, which is of the order
$n^{-\alpha
/(2\alpha+1)}\sqrt{\log n}$.
\end{pf*}

We also wish to study the coverage problem for $L_\infty$-credible
regions when the regularity $\alpha$ is not known. Consider the
empirical Bayes device of the paper under discussion and assume that
the true sequence has a polished tail. The heuristic arguments given
below seem to indicate that the credible region constructed by plugging
in the empirical Bayes estimate of $\alpha$ should have adequate coverage.

Because we deal with various values of $\alpha$ simultaneously, let us
include $\alpha$ in the notation $\Pi_\alpha$ for the prior,
$Z_\alpha$
and $V_\alpha$ for the Gaussian processes introduced in the proof, and
$\epsilon_{n,\alpha}=n^{-\alpha/(2\alpha+1)}\sqrt{\log n}$ for the
sup-norm posterior contraction rate. We observe that $\epsilon
_{n,\alpha
}$ is decreasing in $\alpha$.
By Theorem~5.1 of the paper under discussion, it follows that the
empirical Bayes estimate $\hat{\alpha}$ of $\alpha$ lies, with high
probability, between two deterministic bounds $\underline{\alpha}$ and
$\overline{\alpha}$, and that $\epsilon_{n,\underline{\alpha
}}\asymp
\epsilon_{n,\overline{\alpha}}$.

In the proof of the result on coverage of the credible region, one
needs to lower bound the radius of the credible ball around the
estimate and show that its order is at least as large as the
convergence rate of the point estimator given by the center of the
credible region. When $\hat{\alpha}$ is plugged in, the radius of the
credible region is of the order of the expected value of the supremum
of the Gaussian process $Z_{\hat{\alpha}}$. The randomness of this
process comes from posterior variation conditioned on the sample, and
hence $\hat{\alpha}$ can be considered as a constant. Therefore, as
argued in the proof of the theorem, radius of the credible region is of
the order $\epsilon_{n,\hat{\alpha}}\asymp\epsilon_{n,\underline
{\alpha
}}\asymp
\epsilon_{n,\overline{\alpha}}$.

The sampling error of the Bayes estimator $\hat{f}_\alpha$ using $\Pi
_\alpha$ has two parts---variability around its expectation
$Z_\alpha
$ and its bias. Now for any $t,s\in[0,1]$, $\mathrm{E}|Z_\alpha
(t)-Z_\alpha
(s)|^2$ is decreasing in $\alpha$, so, by Slepian's inequality,
\[
\sup\bigl\{\mathrm{E}\|Z_\alpha\|_\infty\dvtx \underline{\alpha}\le
\alpha\le \overline {\alpha}\bigr\}= \mathrm{E}\|Z_{\underline{\alpha}}\|_\infty
\asymp\epsilon _{n,\underline{\alpha}
}\asymp\epsilon_{n,\overline{\alpha} },
\]
and fixed quantiles of $\|Z_\alpha\|_\infty$ also have the same order as
the expectation of $\|Z_\alpha\|_\infty$ by Borell's inequality. On the
other hand, the bias of $\hat{f}_\alpha$ increases with~$\alpha$, and
hence its maximum is attained at $\overline{\alpha}$ for $\underline
{\alpha}\le\alpha\le\overline{\alpha}$. Note that if $\overline
{\alpha
}$ underestimates the true $\alpha$, then the order of the bias is
$\epsilon_{n,\overline{\alpha} }\asymp\epsilon_{n,\underline
{\alpha}
}$, and so for every $\alpha$ lying in the range $[\underline{\alpha
},\overline{\alpha}]$, the posterior contraction rate would be the same.
Lemma~3.11 seems to indicate that this may be the case. This will
ensure adequate coverage of the empirical Bayes credible set.

Another issue that might be of interest for future investigation is the
handling of unknown variance. In the nonadaptive setting, both
empirical and hierarchical Bayes approaches can fruitfully address the
issue of unknown variance as demonstrated by \citet{Yoo2014} for
nonparametric regression. In the adaptive setting, this is somewhat
unclear, as the empirical Bayes estimate of smoothness and variance
will depend on each other.

It is also natural to ask if the hierarchical Bayes credible sets can
also have adequate coverage in the adaptive setting. This may not have
an affirmative answer, as indicated by \citet{Rivoirard2012}.

Finally, for other curve estimation problems like density estimation or
nonparametric regression, what should be a proper analog of conditions
like self-similarity of polished tail, and how may that help in
establishing coverage? The nonparametric regression problem may be more
tractable than the density estimation, since for the former a basis
expansion approach reduces the function of interest to a sequence of
real-valued parameters which are typically given normal priors as well
and conjugacy holds in the model. Usually it is more convenient to use
a truncated series expansion, but then the sequence of parameters form
a triangular array. It seems that the main challenge will be to
identify a proper analog of a condition on the tail of the sequence in
such a setting.






\printaddresses
\end{document}